\documentclass[twoside,11pt,a4paper,leqno]{article}
\usepackage{amssymb,amsmath,amscd,euscript,verbatim,array}
\usepackage[center,pagestyles]{titlesec}
\usepackage{geometry}
\usepackage{anysize}
\usepackage{fancyhdr}
\usepackage{indentfirst}
\usepackage{graphicx}
\usepackage{color}
\usepackage{ifpdf}
\marginsize{3.5cm}{3.5cm}{2cm}{2cm}

\titleformat{\section}{\centering\normalsize}{\thesection.}{0.5em}{}
\titleformat{\subsection}{\normalsize\bfseries}{\thesubsection.}{0.5em}{}
\titleformat{\subsubsection}{\normalsize\bfseries}{\thesubsubsection.}{0.5em}{}
\newcommand{\N}{\mathbb{N}}
\newcommand{\Z}{\mathbb{Z}}
\newcommand{\R}{\mathbb{R}}

\newtheorem{Theorem}{Theorem}[section]
\newtheorem{Definition}[Theorem]{Definition}
\newtheorem{Lemma}[Theorem]{Lemma}
\newtheorem{Exercise}[Theorem]{Exercise}
\newtheorem{Proposition}[Theorem]{Proposition}
\newtheorem{Remark}[Theorem]{Remark}
\newtheorem{Corollary}[Theorem]{Corollary}

\newcommand{\eps}{\varepsilon}

\newcommand{\T}{\mathbb{T}}
\newcommand{\bthm}{\begin{Theorem}}
\newcommand{\ethm}{\end{Theorem}}
\newcommand{\bpr}{\begin{Proposition}}
\newcommand{\epr}{\end{Proposition}}
\newcommand{\blm}{\begin{Lemma}}
\newcommand{\elm}{\end{Lemma}}
\newcommand{\bex}{\begin{Exercise}}
\newcommand{\eex}{\end{Exercise}}
\newcommand{\be}{\begin{equation}}
\newcommand{\ee}{\end{equation}}
\newcommand{\beal}{\begin{aligned}}
\newcommand{\enal}{\end{aligned}}
\newcommand{\brm}{\begin{Remark}}
\newcommand{\erm}{\end{Remark}}
\newcounter{item}[section]

\newcommand{\Proof}{\textbf{Proof}\hspace{0.3cm}}
\newcommand{\End}{\ensuremath{\hfill{\Box}}\\}
\renewcommand{\title}[1]{\begin{center}\textbf{\large #1}\end{center}}
\renewcommand{\author}[1]{\begin{center}\small #1\end{center}}
\renewcommand{\date}[1]{\begin{center}#1\end{center}}

\makeatletter \@addtoreset{equation}{section}
\makeatother
\begin{document}
\vspace{10pt}
\title{TOTAL DESTRUCTION OF LAGRANGIAN TORI}
\vspace{6pt}
\author{\sc Lin Wang}
\vspace{10pt} \thispagestyle{plain}
{\begingroup\makeatletter
\let\@makefnmark\relax
\makeatother\endgroup}
\begin{quote}
\small {\sc Abstract.} For an integrable Tonelli
Hamiltonian with $d\ (d\geq 2)$ degrees of freedom, we show
that all of the Lagrangian tori can be destroyed by analytic perturbations
which are arbitrarily small in the $C^{d-\delta}$ topology.
\end{quote}
\begin{quote}
\end{quote}
\begin{quote}
\end{quote} \vspace{25pt}

\section{\sc Introduction and main result}
For exact area-preserving twist maps on annulus, it was proved by
Herman in \cite{H2} that invariant circles with given rotation
numbers can be destroyed by $C^{3-\delta}$ arbitrarily small
$C^\infty$ perturbations, where $\delta$ is a small positive
constant. Following the ideas and techniques developed by J.N.Mather
in the series of papers [M1], [M2], [M3] and [M4], a variational
proof of Herman's result was provided in \cite{W1}. A Hamiltonian  is called a Tonelli Hamiltonian if it satisfies positive definiteness, superlinear growth with respect to momentum and completeness of flow. For the case with $d\geq 2$ degrees of
freedom, it was proved in \cite{CW} that for every given rotation
vector, invariant Lagrangian torus with that rotation vector of an
integrable Tonelli Hamiltonian system can be destroyed by an arbitrarily
small $C^\infty$ perturbation in the $C^{2d-\delta}$ topology. In
contrast with it, it was shown that KAM torus with Diophantine
frequency persists under $C^{2d+\delta}$-small perturbations
(\cite{P}). Hence, the above result is  almost optimal.

On the other hand, it was proved that all invariant circles can be
destroyed by $C^1$ arbitrarily small $C^\infty$ perturbations of
the integrable area-preserving twist maps \cite{T}. $C^1$ topology was improved to
be $C^{2-\delta}$ topology by Herman in \cite{H1}. Moreover, he
extended the result to systems with multi-degrees of freedom and
found  that all of the Lagrangian tori of an integrable symplectic twist
map can be destroyed by $C^{d+2-\delta}$ arbitrarily small
$C^\infty$ perturbations of the generating function \cite{H3}. By \cite{H2}, each orbit on an invariant Lagrangian graphs is an action minimizing curve. Based on the minimality of the orbits on the Lagrangian graph, \cite{MMS}  provided a  criterion of non-existence of all Lagrangian tori and applied it to a particular example.
Based on the correspondence between symplectic twist maps and Tonelli Hamiltonian systems (\cite{Go,Mo4}), it shows that all of the Lagrangian tori of
an integrable Tonelli Hamiltonian system with $d\geq 2$ degrees of freedom
can be destroyed by $C^\infty$ perturbations which are arbitrarily
small in the $C^{d+1-\delta}$ topology.

Comparing the results on both sides, it is natural to ask the
following question: \begin{itemize}
\item  if all of the Lagrangian tori can
be destroyed by an arbitrarily small real-analytic perturbation in
the $C^r$ topology, what is the maximum of $r$?
\end{itemize}

In this paper, we prove the following theorem:
\begin{Theorem}\label{Mt}
For an integrable Tonelli Hamiltonian with $d\
(d\geq 2)$ degrees of freedom, all of the Lagrangian tori can be destroyed
by analytic perturbations which are arbitrarily small in the
$C^{d-\delta}$ topology for a small given constant $\delta>0$.
\end{Theorem}

Unfortunately,
 we still don't know
whether the other results are optimal. Some further developments of
KAM theory are needed to verify the optimality. To prove Theorem
\ref{Mt}, we recall some notions on Lagrangian torus.

 In $\text{T}^\ast\T^d$, a submanifold $\mathcal {T}^d$ is called
Lagrangian torus if it is diffeomorphic to the torus $\T^d$ and the
symplectic form (non-degenerate closed 2-form) vanishes on it. An example of Lagrangian torus is
the KAM torus.
\begin{Definition}\label{dd}$\bar{\mathcal
{T}}^d$ is called a $d$ dimensional KAM torus if
\begin{itemize}
\item $\bar{\mathcal {T}}^d$ is a Lipschitz graph over $\T^d$;
\item $\bar{\mathcal {T}}^d$ is invariant under the Hamiltonian flow
 $\Phi_t^H$ generated by the Hamiltonian function $H$;
\item there exists a diffeomorphism
$\phi:\ \T^d\rightarrow \bar{\mathcal {T}}^d$ such that
$\phi^{-1}\circ\Phi_H^t\circ\phi=R_\omega^t$ for any $t\in \R$,
where $R_\omega^t:\ x\rightarrow x+\omega t$ and $\omega$ is called
the rotation vector of $\bar{\mathcal {T}}^d$.
\end{itemize}
\end{Definition}
Generally, the rotation vector of the Lagrangian torus
 is not well defined if the Lagrangian torus contains several invariant sets with different rotation
vectors. In this paper, we are concerned with the Lagrangian torus as
follow:
\begin{Definition}\label{dd1}
$\mathcal {T}^d$ is called a $d$ dimensional Lagrangian torus if
\begin{itemize}
\item $\mathcal {T}^d$ is a Lagrangian graph;
\item $\mathcal {T}^d$ is invariant for the Hamiltonian
flow $\Phi_H^t$ generated by $H$.
\end{itemize}
\end{Definition}

 It is still open whether all the invariant Lagrangian tori of Tonelli Hamiltonians are graphs or not. Some results have been obtained by adding topological or variational conditions (see \cite{Arna}, \cite{BP} and \cite{CR04} for instance). \cite{Arna} proved that for Tonelli Hamiltonians, the invariant Lagrangian submanifolds that are isotopic to zero section are graphs. It is easy to see that ``isotopic to zero section" is necessary if we consider the elliptic island for the Hamiltonian time one map of rigid pendulum. However, the orbits on the elliptic island is not minimal in the sense of variation. Under the condition that every orbit is an action minimizing curve, \cite{CR04} proved that for  Tonelli Hamiltonians with 2 degrees of freedom, the invariant Lagrangian tori are graphs.  By \cite{H2}, each orbit on an invariant Lagrangian graph is an action minimizing curve. Hence, it is shown that the graph property of a Lagrangian torus is equivalent to the minimality of the orbits on the torus for the case with two degrees of freedom. Whereas, it is still open to verify the equivalence for the case with multi-degrees of freedom.

In the following sections, we consider the destruction of all of the
Lagrangian tori of symplectic twist maps.
\begin{Definition}\label{dd3}
A map $f:\,(x,y)\rightarrow (x',y')$ of $\T^d\times\R^d$ is called a symplectic twist map if
\begin{itemize}
\item $f$ is a diffeomorphism isotopic to the identity;
\item $f$ preserves the symplectic form;
\item $\frac{\partial x'}{\partial y}$ is uniformly positive definite and bounded.
\end{itemize}
\end{Definition}
In particular, the last one is called twist condition. An necessary condition for existence of any invariant Lagrangian graph is that there exists a generating function $h:\, \T^d\times\R^d\rightarrow \R$ such that
\[y'=\partial_2h(x,x')\qquad y=-\partial_1h(x,x')\]
where $\partial_1$, $\partial_2$ denote derivatives with respect to the first and second arguments respectively. We will restrict attention to maps with periodic generating functions, also known as exact symplectic maps.

Based on the
correspondence between symplectic twist maps and Hamiltonian
systems, it can be achieved to destroy all of the Lagrangian tori of
Hamiltonian systems.

\section{\sc A toy model} To show the basic ideas, we are beginning
with a toy model whose generating function is as follow:
\begin{equation}\label{hall}
h_n(x,x')=h_0(x,x')-\frac{5}{4n^2}\sin(nx')-\frac{1}{16 n^2}\cos(2
nx'),
\end{equation}
where $h_0(x,x')=\frac{1}{2}(x-x')^2$. Let $f_n(x,y)=(x',y')$ be the
exact area-preserving twist map generated by (\ref{hall}), then
\begin{equation*}
\begin{cases}
y=-\partial_1 h_n(x,x')=x'-x,\\
y'=\partial_2 h_n(x,x')=x'-x-\frac{5}{4 n}\cos( nx')+\frac{1}{8
n}\sin(2 nx').
\end{cases}
\end{equation*}
We set $\phi_n(x)=-\frac{5}{4 n}\cos( nx)+\frac{1}{8 n}\sin(2 nx)$,
then
\begin{equation}
f_n(x,y)=(x+y,y+\phi_n(x+y)).
\end{equation}
In \cite{H1}, Herman found a criterion of total destruction of
invariant circles. By Birkhoff graph theorem (see \cite{H2}), if
$f_n$ admits an invariant circle, then the invariant circle is a
Lipschitz graph. We denote the graph by $\psi_n$, then it follows
from \cite{H3} that
\begin{equation*}
\psi_n\circ g_n-\psi_n=\phi_n\circ g_n,
\end{equation*}
where $g_n=\text{Id}+\psi_n$. This is equivalent to
\begin{equation}\label{gg}
\frac{1}{2}(g_n+g_n^{-1})=\text{Id}+\frac{1}{2}\phi_n.
\end{equation}
Let $\mathfrak{D}_n$ be the set of differentiate points of $g_n$,
then $\mathfrak{D}_n$ has full Lebesgue measure on $\R$ since $g_n$
is a Lipschitz function. For $x\in \mathfrak{D}_n$, we differentiae
(\ref{gg}),
\begin{equation*}
\frac{1}{2}(Dg_n(x)+(Dg_n)^{-1}(g_n^{-1}(x)))=1+\frac{1}{2}D\phi_n(x).
\end{equation*}
Let $G_n=||Dg_n||_{L^\infty}$. It is easy to see that for $\eps>0$,
there exists $\tilde{x}\in \mathfrak{D}_n$ such that
$Dg_n(\tilde{x})\geq G_n-\eps$. Let $M_n=\max D\phi_n$, we have
\begin{equation*}
\frac{1}{2}\left(G_n+\frac{1}{G_n}-\eps\right)\leq 1+\frac{1}{2}M_n.
\end{equation*}
Since $\eps>0$ is arbitrarily small, then
\begin{equation*}
\frac{1}{2}\left(G_n+\frac{1}{G_n}\right)\leq 1+\frac{1}{2}M_n.
\end{equation*}
Hence,
\begin{equation}\label{mm}
G_n\leq
1+\frac{1}{2}M_n+\left(M_n+\frac{1}{4}{M_n}^2\right)^{\frac{1}{2}}.
\end{equation}
Obviously, for $x\in \mathfrak{D}_n$, we have
\begin{equation*}
\frac{1}{G_n}\leq 1+\frac{1}{2}D\phi_n(x).
\end{equation*}
Let $m_n=\min D\phi_n$, then we have
\begin{equation*}
\frac{1}{G_n}\leq 1+\frac{1}{2}m_n,
\end{equation*}
which together with (\ref{mm}) implies that
\begin{equation*}
\frac{1}{1+\frac{1}{2}m_n}\leq
1+\frac{1}{2}M_n+\left(M_n+\frac{1}{4}{M_n}^2\right)^{\frac{1}{2}}.
\end{equation*}
Therefore, it is sufficient for total destruction of invariant
circles to construct $\phi_n(x)$ such that
\begin{equation}\label{cri}
\frac{1}{1+\frac{1}{2}\min D\phi_n}> 1+\frac{1}{2}\max
D\phi_n+\left(\max D\phi_n+\frac{1}{4}(\max
D\phi_n)^2\right)^{\frac{1}{2}}.
\end{equation}
In our construction,
\begin{equation}
D\phi_n(x)=\frac{5}{4}\sin(nx)+\frac{1}{4}\cos(2 nx).
\end{equation}
A simple calculation implies
\begin{equation*}
\left\{\begin{array}{ll}
\hspace{-0.4em}\min D\phi_n(x)=-\frac{3}{2},&\text{attained at}\ x=\frac{3}{2n}+\frac{2\pi k}{n},\\
\hspace{-0.4em}\max D\phi_n(x)=1,&\text{attained at}\ x=\frac{\pi}{2n}+\frac{2\pi k}{n},\\
\end{array}\right.
\end{equation*}
where $k\in\Z$. Hence, (\ref{cri}) holds. Moreover, the exact
area-preserving twist map generated by (\ref{hall}) admits no
invariant circles.

By interpolation inequality (\cite{H1}), for a small positive
constant $\delta$, we have
\[||\phi_n||_{C^{1-\delta}}\leq 2||\phi_n||_{C^0}^\delta||D\phi_n||_{C^0}^{1-\delta}.\]
From the construction of $\phi_n$, it follows that
$||\phi_n||_{C^0}\rightarrow 0$, as $n\rightarrow \infty$ and
$||D\phi_n||_{C^0}$ is bounded. Hence,
\begin{equation*}
||\phi_n||_{C^{1-\delta}}\rightarrow 0\quad\text{as}\quad
n\rightarrow\infty,
\end{equation*}
which implies that
\begin{equation*}
||h_n-h_0||_{C^{2-\delta}}\rightarrow 0\quad\text{as}\quad
n\rightarrow\infty.
\end{equation*}
\section{\sc $C^\infty$ destruction of all of the Lagrangian tori}
 In \cite{H3}, Herman extended
the criterion (\ref{cri}) to multi-degrees of freedom. More
precisely, for exact symplectic twist map on $\text{T}^*\T^d$, whose generating function is
\begin{equation}
h(x,x')=\frac{1}{2}(x-x')^2+\Psi(x'),
\end{equation}
where $\Psi\in C^r(\T^d,\R)$, $r\geq 2$. Correspondingly,  the exact symplectic twist map has the following form
\begin{equation}
f(x,y)=(x+y,y+d\Psi(x+y)),
\end{equation}
where
\[d\Psi=\left(\frac{\partial\Psi}{\partial x_1},\cdots,\frac{\partial\Psi}{\partial x_d}\right).\]
Let $E(x)$ be the derivative matrix of $d\Psi$ and
$T(x)=\frac{1}{d}\text{tr}E(x)$, where $\text{tr}E(x)$ denotes the trace
of $E(x)$. From a similar argument as the deduction of (\ref{cri}),
it follows that it is sufficient for total destruction of of the Lagrangian
tori to construct $T(x)$ such that
\begin{equation}\label{cri_multi}
\frac{1}{1+\frac{1}{2}\min T(x)}> 1+\frac{1}{2}\max T(x)+\left(\max
T(x)+\frac{1}{4}(\max T(x))^2\right)^{\frac{1}{2}}.
\end{equation}
Moreover, for $T(x)\rightarrow 0$, (\ref{cri_multi}) implies
\begin{equation}\label{scri}
-\frac{1}{2}\min T(x)>\sqrt{\max T(x)}+O(\max T(x)).
\end{equation}

\begin{Remark}
The integrable part $\frac{1}{2}(x-x')^2$ can be easily generalized to the form $\frac{1}{2}(x-x')^tA(x-x')$, where $(\cdot)^t$ denotes the transpose of $(\cdot)$ and $A$ denotes a symmetric positive definite matrix. By a similar calculation, (\ref{cri_multi})  still holds true if $T(x)=\frac{1}{d}\text{tr}E(x)$ is replaced by  $T(x)=\frac{1}{d}\text{tr}\left(A^{-1}E(x)\right)$. The choice of the form $\frac{1}{2}(x-x')^2$ could set us
free from a tedious calculation to see the crucial mechanism of the
problem.
\end{Remark}

Herman constructed a sequence $\{\Psi_n\}_{n\in\N}$ that satisfies
(\ref{scri}). It is easy to see $T_n(x)=\frac{1}{d}\Delta\Psi_n$
where $\Delta$ denotes the Laplacian. Since $T_n(x)$ is
$2\pi$-periodic, it is enough to construct it on $[-\pi,\pi]^d$.
More precisely,
\begin{equation*}
T_n(x)=\left\{\begin{array}{ll}
\hspace{-0.4em}T_n^+(x),&x\in [0,\pi]^d,\\
\hspace{-0.4em}-T_n^-(x),&x\in [-\pi,0]^d,\\
\hspace{-0.4em}0,&\text{others}.\\
\end{array}\right.
\end{equation*}
where $T_n(x)$ is $C^\infty$ function, $T_n^+(x)$ and $T_n^-(x)$
have the following forms respectively.

$T_n^+(x)$ satisfies:
\begin{equation*}
\left\{\begin{array}{l}
\hspace{-0.4em}\text{supp}\,T_n^+(x)\subset [0,\pi]^d,\\
\hspace{-0.4em}\max T_n^+(x)=\frac{1}{9n},\\
\end{array}\right.
\end{equation*}

$T_n^-(x)$ satisfies:
\begin{equation*}
\left\{\begin{array}{l}
\hspace{-0.4em}\text{supp}\,T_n^-(x)=B_{R_n}(x_0),\\
\hspace{-0.4em}\max T_n^-(x)=\frac{1}{\sqrt{n}},\\
\hspace{-0.4em}R_n\sim
\left(\frac{1}{\sqrt{n}}\right)^\frac{1}{d},\\
\hspace{-0.4em}x_0=\left(-\frac{\pi}{2},\cdots,-\frac{\pi}{2}\right),\\
\end{array}\right.
\end{equation*}
where $f\sim g$ means that $\frac{1}{C}g<f<Cg$ holds for a constant
$C>1$. Hence, we obtain a sequence of $\{T_n(x)\}_{n\in \N}$ with
bounded $C^d$ norms and satisfying
\[\int_{\T^d}T_n(x)dx=0.\]
From interpolation inequality, it follows that $T_n(x)\rightarrow 0$
as $n\rightarrow \infty$ in the $C^{d-\delta}$ topology for any
$\delta>0$.

Let $\Psi_n$ be the unique function in $C^\infty(\T^d,\R)$ such that
\[\int_{\T^d}\Psi_n(x)dx=0\quad\text{and}\quad\frac{1}{d}\Delta\Psi_n(x)=T_n(x).\]
By Schauder estimates one knows that for any $\delta>0$,
$\Psi_n(x)\rightarrow 0$ as $n\rightarrow\infty$ in the
$C^{d+2-\delta}$ topology. From the construction of $T_n(x)$, it is
easy to see that (\ref{scri}) is verified.

Above all, we have the following theorem
\begin{Theorem}
All of the Lagrangian tori of an integrable symplectic twist map with
$d\geq 1$ degrees of freedom can be destroyed by $C^\infty$
perturbations of the generating function and the perturbations are
arbitrarily small in the $C^{d+2-\delta}$  topology for a small
given constant $\delta>0$.
\end{Theorem}
Based on the correspondence between symplectic twist maps and
Hamiltonian systems, we have the following corollary.
\begin{Corollary}
All of the Lagrangian tori of an integrable Tonelli Hamiltonian system with $d\geq
2$ degrees of freedom can be destroyed by $C^\infty$ perturbations
which are arbitrarily small in the $C^{d+1-\delta}$ topology  for a
small given constant $\delta>0$.
\end{Corollary}
\section{\sc An approximation lemma}
In this section, we will prove a lemma on $C^\infty$ functions
approximated by trigonometric polynomials. First of all, we need
some notations. Define
\[C^\infty_{2\pi}(\R^d,\R):=\left\{f:\R^d\rightarrow\R|f\in C^\infty(\R^d,\R)\ \text{and}\ 2\pi-\text{periodic in}\ x_1,\ldots,x_d\right\}.\]
 Let $f(x)\in C^\infty_{2\pi}(\R^d,\R)$. The $m$-th
 Fej\'{e}r-polynomial of $f$ with respect to $x_j$ is given by
 \begin{equation}\label{B1}
F_m^{[j]}(f)(x):=\frac{1}{m\pi}\int_{-\pi/2}^{\pi/2}f(x+2te_j)\left(\frac{\sin(mt)}{\sin
t}\right)^2dt,
 \end{equation}
where $x\in\R^d$, $m\in\N$, $j\in\{1,\ldots,d\}$ and $e_j$ is the
$j$-th vector of the canonical basis of $\R^d$. $F_m^{[j]}(f)(x)$ is
a trigonometric polynomial in $x_j$ of degree at most $m-1$. By
\cite{Z},
\[\frac{1}{m\pi}\int_{-\pi/2}^{\pi/2}\left(\frac{\sin(mt)}{\sin
t}\right)^2dt=1,\] hence, from (\ref{B1}), we have
\[||F_m^{[j]}(f)||_{C^0}\leq ||f||_{C^0}.\]
We denote
\[P_m^{[j]}(f):=2F_{2m}^{[j]}(f)-F_m^{[j]}(f).\] It is easy to see
that $P_m^{[j]}(f)$ is a trigonometric polynomial in $x_j$ of degree
at most $2m-1$. Moreover,
 \begin{equation}\label{B2}
||P_m^{[j]}(f)||_{C^0}\leq 3||f||_{C^0},
 \end{equation}
 \begin{equation}\label{B3}
P_m^{[j]}(af+bg)=aP_m^{[j]}(f)+bP_m^{[j]}(g),
 \end{equation}
where $a,b\in\R$ and $f,g\in C^\infty_{2\pi}(\R^d,\R)$. For
$k\in\{1,\ldots,d\}$, $j_1,\ldots,j_k\in\{1,\ldots,d\}$ with
$j_p\neq j_q$ for $p\neq q$. Let $m_1,\ldots,m_k\in\N$ and $f\in
C^\infty_{2\pi}(\R^d,\R)$, we define
 \begin{equation}\label{B4}
P_{m_1,\ldots,m_k}^{[j_1,\ldots,j_k]}(f):=P_{m_1}^{[j_1]}\left(P_{m_2}^{[j_2]}\left(\cdots\left(P_{m_k}^{[j_k]}(f)\right)\cdots\right)\right).
 \end{equation}
It is easy to see that for all $l\in \{1,\ldots,k\}$,
$P_{m_1,\ldots,m_k}^{[j_1,\ldots,j_k]}(f)$ are trigonometric
polynomials in $x_{j_l}$ of degree at most $2m_l-1$, also known as generalized de la Vall\'{e}e Poussin polynomial. We have the
following lemma.
\begin{Lemma}\label{appp}
Let $f\in C^\infty_{2\pi}(\R^d,\R)$, $r_1,\ldots,r_d\in\N$,
$m_1,\ldots,m_d\in\N$, then we have
\begin{equation}\label{B8}
\|f-P_{m_1,\ldots,m_d}^{[1,\ldots,d]}(f)\|_{C^0}\leq
C_d\sum_{j=1}^{d}\frac{1}{{m_j}^{r_j}}\left\|\frac{\partial
^{r_j}f}{\partial {x_j}^{r_j}}\right\|_{C^0},
\end{equation}
where $C_d$ is a constant only depending on $d$.
\end{Lemma}
Lemma \ref{appp} is a direct corollary of Theorem 2.12 of \cite{A}. For the sake of completeness, we decide to provide the proof.

\Proof We will prove Lemma \ref{appp} by induction.
The case $d=1$ is covered by Jackson's approximation theorem. More
precisely, for $f\in C^\infty_{2\pi}(\R,\R)$, $m,r\in\N$, we have
\begin{equation}\label{B5}
\|f-P_m^{[1]}(f)\|_{C^0}\leq
C_1\frac{1}{{m}^{r}}\left\|\frac{\partial ^{r}f}{\partial
{x}^{r}}\right\|_{C^0}.
\end{equation}
Let the assertion be true for some $d\in\N$. We verify it for $d+1$.
Consider the functions $f(x_1,\cdot)$ with $x_1$ as a real
parameter. Then by the assertion for $d$, we have
\[\|f(x_1,\cdot)-P_{m_2,\ldots,m_{d+1}}^{[2,\ldots,d+1]}(f)(x_1,\cdot)\|_{C^0}\leq
C_d\sum_{j=2}^{d+1}\frac{1}{{m_j}^{r_j}}\left\|\frac{\partial
^{r_j}f}{\partial {x_j}^{r_j}}\right\|_{C^0},\] hence,
\begin{equation}\label{B6}
\|f-P_{m_2,\ldots,m_{d+1}}^{[2,\ldots,d+1]}(f)\|_{C^0}\leq
C_d\sum_{j=2}^{d+1}\frac{1}{{m_j}^{r_j}}\left\|\frac{\partial
^{r_j}f}{\partial {x_j}^{r_j}}\right\|_{C^0}.
\end{equation}
Let $\hat{x}_j\in\R^d$ denote the vector $x\in\R^{d+1}$ without its
$j$-th entry. For the functions $f(\cdot,\hat{x}_1)$, from
(\ref{B5}), it follows that
\[\|f(\cdot,\hat{x}_1)-P_{m_1}^{[1]}(f)(\cdot,\hat{x}_1)\|_{C^0}\leq C_1\frac{1}{{m_1}^{r_1}}\left\|\frac{\partial ^{r_1}f}{\partial
{x_1}^{r_1}}\right\|_{C^0},\]hence,
\begin{equation}\label{B7}
\|f-P_{m_1}^{[1]}(f)\|_{C^0}\leq
C_1\frac{1}{{m_1}^{r_1}}\left\|\frac{\partial ^{r_1}f}{\partial
{x_1}^{r_1}}\right\|_{C^0}.
\end{equation}
By (\ref{B2}), (\ref{B3})),(\ref{B4}) and (\ref{B6}), we have
\begin{align*}
\left\|P_{m_1}^{[1]}(f)-P_{m_1,\ldots,m_{d+1}}^{[1,\ldots,d+1]}(f)\right\|_{C^0}&=\left\|P_{m_1}^{[1]}(f)-P_{m_1}^{[1]}\left(P_{m_2,\ldots,m_{d+1}}^{
[2,\ldots,d+1]}(f)\right)\right\|_{C^0},\\
&=\left\|P_{m_1}^{[1]}\left(f-P_{m_1}^{[1]}P_{m_2,\ldots,m_{d+1}}^{
[2,\ldots,d+1]}(f)\right)\right\|_{C^0},\\
&\leq 3\left\|f-P_{m_1}^{[1]}P_{m_2,\ldots,m_{d+1}}^{
[2,\ldots,d+1]}(f)\right\|_{C^0},\\
&\leq 3C_d\sum_{j=2}^{d+1}\frac{1}{{m_j}^{r_j}}\left\|\frac{\partial
^{r_j}f}{\partial {x_j}^{r_j}}\right\|_{C^0},
\end{align*}
which together with (\ref{B7}) implies that
\begin{align*}
\left\|f-P_{m_1,\ldots,m_{d+1}}^{[1,\ldots,d+1]}(f)\right\|_{C^0}&\leq \left\|f-P_{m_1}^{[1]}(f)\right\|_{C^0}+\left\|P_{m_1}^{[1]}(f)-P_{m_1,\ldots,m_{d+1}}^{[1,\ldots,d+1]}(f)\right\|_{C^0},\\
&=\left\|f-P_{m_1}^{[1]}(f)\right\|_{C^0}+\left\|P_{m_1}^{[1]}(f)-P_{m_1}^{[1]}\left(P_{m_2,\ldots,m_{d+1}}^{
[2,\ldots,d+1]}(f)\right)\right\|_{C^0},\\
&\leq C_1\frac{1}{{m_1}^{r_1}}\left\|\frac{\partial
^{r_1}f}{\partial
{x_1}^{r_1}}\right\|_{C^0}+3C_d\sum_{j=2}^{d+1}\frac{1}{{m_j}^{r_j}}\left\|\frac{\partial
^{r_j}f}{\partial {x_j}^{r_j}}\right\|_{C^0},\\
&\leq
C_{d+1}\sum_{j=1}^{d+1}\frac{1}{{m_j}^{r_j}}\left\|\frac{\partial
^{r_j}f}{\partial {x_j}^{r_j}}\right\|_{C^0}.
\end{align*}
This finishes the proof of Lemma \ref{appp}. \End

Obviously, there exist $m_{\bar{j}}, r_{\bar{j}}$ such that
\[\frac{1}{{m_{\bar{j}}}^{r_{\bar{j}}}}\left\|\frac{\partial
^{r_{\bar{j}}}f}{\partial
{x_{\bar{j}}}^{r_{\bar{j}}}}\right\|_{C^0}=\max_{1\leq j \leq
d}\left\{\frac{1}{{m_j}^{r_j}}\left\|\frac{\partial
^{r_j}f}{\partial {x_j}^{r_j}}\right\|_{C^0}\right\}.\] Hence, we
have
\begin{align*}
\|f-P_{m_1,\ldots,m_d}^{[1,\ldots,d]}(f)\|_{C^0}&\leq
dC_d\frac{1}{{m_{\bar{j}}}^{r_{\bar{j}}}}\left\|\frac{\partial
^{r_{\bar{j}}}f}{\partial
{x_{\bar{j}}}^{r_{\bar{j}}}}\right\|_{C^0},\\
&\leq
C'_d\frac{1}{{m_{\bar{j}}}^{r_{\bar{j}}}}\|f\|_{C^{r_{\bar{j}}}}.
\end{align*}
For the simplicity of notations, we denote
\[p_N(x)=P_{m_1,\ldots,m_d}^{[1,\ldots,d]}(f)(x),\]
where $x=(x_1,\ldots,x_d)$ and $N=2m_{\bar{j}}-1$. Moreover, we
denote $k:=r_{\bar{j}}$, then
\begin{equation}\label{B9}
\|f(x)-p_N(x)\|_{C^0}\leq A_{dk}N^{-k}\|f(x)\|_{C^k},
\end{equation}
where $A_{dk}$ is a constant depending on $d$ and $k$.

\section{\sc $C^\omega$ destruction of all of the Lagrangian tori}
Similar to Herman's construction,  we consider $C^\infty$ function
$\tilde{T}_n(x)$ as follow:
\begin{equation*}
\tilde{T}_n(x)=\left\{\begin{array}{ll}
\hspace{-0.4em}\tilde{T}_n^+(x),&x\in [0,\pi]^d,\\
\hspace{-0.4em}-\tilde{T}_n^-(x),&x\in [-\pi,0]^d,\\
\hspace{-0.4em}0,&\text{others}.\\
\end{array}\right.
\end{equation*}

$\tilde{T}_n^+(x)$ satisfies:
\begin{equation*}
\left\{\begin{array}{l}
\hspace{-0.4em}\text{supp}\,\tilde{T}_n^+(x)\subset [0,\pi]^d,\\
\hspace{-0.4em}\max \tilde{T}_n^+(x)= 1,\\
\end{array}\right.
\end{equation*}

$\tilde{T}_n^-(x)$ satisfies:
\begin{equation*}
\left\{\begin{array}{l}
\hspace{-0.4em}\text{supp}\,\tilde{T}_n^-(x)=B_{R_n}(x_0),\\
\hspace{-0.4em}\max \tilde{T}_n^-(x)=n,\\
\hspace{-0.4em}R_n\sim
\left(\frac{1}{n}\right)^\frac{1}{d},\\
\hspace{-0.4em}x_0=\left(-\frac{\pi}{2},\cdots,-\frac{\pi}{2}\right).\\
\end{array}\right.
\end{equation*}

Moreover, we require $\int_{\T^d}\tilde{T}_n(x)dx=0$.
By Lemma \ref{appp}, there exists a trigonometric polynomial
$p_N(x_1,\cdots,x_d)$ in $x_l$ $(1\leq l \leq d)$ of degree at most
$N$ such that
\begin{equation}\label{Moti}
\left\|\tilde{T}_n(x_1,\cdots,x_d)-p_N(x_1,\cdots,x_d)\right\|_{C^0}\leq A_{dk}N^{-k}\left\|\tilde{T}_n(x_1,\cdots,x_d)\right\|_{C^k}.
\end{equation}
By the construction of $\tilde{T}_n$, we have
\begin{equation}\label{Tck}
||\tilde{T}_n(x_1,\cdots,x_d)||_{C^k}\sim n^{\frac{k}{d}+1}.
\end{equation}
Then, choosing $N$ large enough such that
\begin{equation}\label{N}
A_{dk}N^{-k}||\tilde{T}_n(x_1,\cdots,x_d)||_{C^k}<\sigma\ll 1,
\end{equation}
where $\sigma$ is a small enough positive constant. Hence, we have
\begin{equation}
\left\{\begin{array}{ll}\hspace{-0.4em}\max p_N(x)\sim n,&\text{attained on}\ B_{R_n}(x_0),\\
\hspace{-0.4em}\max p_N(x)\sim 1,&\text{on}\ [-\pi,\pi]^d\backslash\ B_{R_n}(x_0).\\
\end{array}\right.
\end{equation}
By (\ref{Tck}) and (\ref{N}), we have
\begin{equation}\label{nn}
N>\left(\frac{A_{dk}}{\sigma}\right)^{\frac{1}{k}}n^{\frac{1}{d}+\frac{1}{k}}\geq
Cn^{\frac{1}{d}+\frac{1}{k}},
\end{equation}
where $C$ is a constant independent of $n$. Since $\int_{\T^d}\tilde{T}_n(x)dx=0$, it follows from (\ref{Moti}) that $\int_{\T^d}p_N(x)dx=0$. We consider the
normalized trigonometric polynomial
\begin{equation}
\tilde{p}_N(x)=\frac{1}{n^{1-\eps}}\frac{p_N(x)}{\max |p_N(x)|},
\end{equation}where $x=(x_1,\ldots,x_d)$. It is easy to see that $\int_{\T^d}\tilde{p}_N(x)dx=0$ and
\begin{equation}\label{mxmi}
\left\{\begin{array}{l}
\hspace{-0.4em}\max \tilde{p}_N(x)\sim \frac{1}{n^{2-\eps}},\\
\hspace{-0.4em}\min \tilde{p}_N(x)\sim -\frac{1}{n^{1-\eps}}.\\
\end{array}\right.
\end{equation}
It follows from (\ref{mxmi}) that
\begin{equation*}
-\frac{1}{2}\min \tilde{p}_N(x)\sim \frac{1}{n^{1-\eps}},
\end{equation*}
\begin{equation*}
\sqrt{\max \tilde{p}_N(x)}+O(\max \tilde{p}_N(x))\sim
\frac{1}{n^{1-\frac{\eps}{2}}}.
\end{equation*}
Hence, for $n$ large enough, we have
\begin{equation}\label{ccrri}
-\frac{1}{2}\min \tilde{p}_N(x)>\sqrt{\max \tilde{p}_N(x)}+O(\max
\tilde{p}_N(x)).
\end{equation}

 Next, we estimate $||\tilde{p}_N(x)||_{C^r}$. By a simple
calculation, we have
\begin{equation}
||p_N(x)||_{C^r}\leq CnN^{r+1}.
\end{equation}
Then,
\begin{align*}
||\tilde{p}_N(x)||_{C^r}&\leq\frac{1}{n^{1-\eps}}\frac{1}{\max
|p_N(x)|} ||p_N(x)||_{C^r},\\
&\leq \frac{1}{n^{2-\eps}}\cdot CnN^{r+1},\\
&=C\frac{1}{n^{1-\eps}}N^{r+1}.
\end{align*}
To achieve $n^{-(1-\eps)}N^{r+1}\rightarrow 0$ as
$n\rightarrow\infty$, it suffices to make
\[\frac{1}{n^{1-\eps}}N^{r+1}\leq \frac{1}{n^\eps}.\] Hence, we have
\[r\leq \log_{N}n^{1-2\eps}-1.\]
From (\ref{nn}), we take
\[N\sim n^{\frac{1}{d}+\frac{1}{k}},\]then, we have
\begin{equation}
\max_N r\leq \left(\frac{1}{d}+\frac{1}{k}\right)^{-1}(1-2\eps)-1.
\end{equation}
Since $k$ can be made large enough, we take $k=\frac{d}{2\eps}$. Let
\[\delta=\frac{4\eps}{1+2\eps}d,\]then, we have
\begin{equation}\label{maxr}
\max_N r\leq d-1-\delta.
\end{equation}
It follows from (\ref{maxr}) that $\tilde{p}_N(x)\rightarrow 0$ as
$n\rightarrow\infty$ in the $C^{d-1-\delta}$ topology for any
$\delta>0$.

Since $\int_{\T^d}\tilde{p}_N(x)dx=0$, then there exists the unique function $\tilde{\Psi}_n\in C^\omega(\T^d,\R)$
such that
\[\int_{\T^d}\tilde{\Psi}_n(x)dx=0\quad\text{and}\quad\frac{1}{d}\Delta\tilde{\Psi}_n(x)=\tilde{p}_N(x).\]
By Schauder estimates one knows that for any $\delta>0$,
$\tilde{\Psi}_n(x)\rightarrow 0$ as $n\rightarrow\infty$ in the
$C^{d+1-\delta}$ topology. According to (\ref{ccrri}),  we have that
the symplectic twist maps generated by the generating function
$\tilde{h}_n(x,x')=\frac{1}{2}(x-x')^2+\tilde{\Psi}_n(x')$ do not
admit any Lagrangian tori for $n$ large enough.

So far, we prove the following theorem
\begin{Theorem}
All of the Lagrangian tori of an integrable symplectic twist map with
$d\geq 2$ degrees of freedom can be destroyed by $C^\omega$
perturbations of the generating function and the perturbations are
arbitrarily small in the $C^{d+1-\delta}$ topology  for a small
given constant $\delta>0$.
\end{Theorem}
Based on the correspondence between symplectic twist maps and
Hamiltonian systems, together with the toy model corresponding to
the case with $d=1$, we finish the proof of Theorem \ref{Mt}.

 \vspace{2ex}
\noindent\textbf{Acknowledgement}  The author sincerely
thanks the referees for their careful reading of the manuscript and
invaluable comments which were very helpful in improving this paper. The author also would like to thank
Prof. C.-Q. Cheng and Dr. L. Jin for many helpful discussions. This
work is under the support of 
the NNSF of
China (Grant 11171146).

{\sc School of Mathematical Sciences, Fudan University,
Shanghai 200433,
China.}

 {\it E-mail address:} \texttt{linwang.math@gmail.com}

\end{document}